\documentclass[a4paper,11pt,oneside]{article}
\usepackage[english]{babel}
\usepackage{amsmath,amsfonts,latexsym,amssymb,amsthm}
\begin{document}
\begin{center}{\Large\bf Complete classification of torsion of elliptic curves over quadratic cyclotomic fields}\\
\vspace{1cm}
Filip Najman\\
\vspace{0.5 cm}
\emph{\small Department of Mathematics, University of Zagreb,\\ Bijeni\v cka cesta 30, 10000 Zagreb, Croatia }
\end{center}
\vspace{1cm}
\begin{abstract} In a previous paper (\cite{fn}), the author examined the possible torsions of an elliptic curve over the quadratic fields $\mathbb Q(i)$ and $\mathbb Q(\sqrt{-3})$.
Although all the possible torsions were found if the elliptic curve has rational coefficients, we were unable to eliminate some possibilities for the torsion if the elliptic curve has coefficients that are not rational. In this note, by finding all the points of two hyperelliptic curves over $\mathbb Q(i)$ and $\mathbb Q(\sqrt{-3})$, we solve this problem completely and thus obtain a classification of all possible torsions of elliptic curves over $\mathbb Q(i)$ and $\mathbb Q(\sqrt{-3})$.
\end{abstract}
\footnotetext{Mathematics subject classification (2000) 11G05, 11G30, 14H40}
\vspace{1cm}
\textbf{1 Introduction}\\
\\
For an elliptic curve $E$ over a number field $K$, it is well known, by the Mordell-Weil theorem, that the set $E(K)$ of $K$-rational points on $E$ is a finitely generated Abelian group. The group $E(K)$ is isomorphic to $T\oplus\mathbb Z^r$, where $r$ is a non-negative integer and $T$ is the torsion subgroup. When $K=\mathbb Q$, by Mazur's Theorem, the torsion subgroup is either cyclic of order $m$, where $1 \leq m \leq 10$ or $m=12$, or of the form $\mathbb Z_2 \oplus \mathbb Z_{2m}$, where $1 \leq m \leq 4$.\\
If $K$ is a quadratic field, then the following theorem classifies the possible torsions.
\newtheorem*{tmx}{Theorem(Kamienny, \cite{kam}, Kenku and Momose, \cite{km})} 
\begin{tmx}
Let $K$ be a quadratic field and $E$ an elliptic curve over $K$. Then the torsion subgroup $E(K)_{tors}$ of $E(K)$ is isomorphic to one of the following $26$ groups:
$$\mathbb Z_m, \text{ for } 1 \leq m\leq 18,\ m\neq 17,$$
$$\mathbb Z_2 \oplus \mathbb Z_{2m}, \text{ for } 1 \leq m\leq 6,$$
$$\mathbb Z_3 \oplus \mathbb Z_{3m}, \text{ for }  m=1,2$$
$$\mathbb Z_4 \oplus \mathbb Z_{4}.$$
\end{tmx}
Moreover, the only quadratic field over which torsion $\mathbb Z_4 \oplus \mathbb Z_4$ appears is $\mathbb Q(i)$ and the only quadratic field over which torsions $\mathbb Z_3 \oplus \mathbb Z_3$ and $\mathbb Z_3 \oplus \mathbb Z_6$ appear is $\mathbb Q(\sqrt{-3})$.\\
In \cite{jkp}, Theorem 3.5, it is proved that if we let the quadratic fields vary, then all of the 26 torsion subgroups appear infinitely often.\\
Unfortunately, if we fix a quadratic field, this theorem does not tell us which of the 26 listed groups actually appear as torsion subgroups. In \cite{fn}, we took this approach, fixing the fields $\mathbb Q(i)$ and $\mathbb Q(\sqrt{-3})$. These fields are somewhat special among quadratic fields, as they are the only cyclotomic quadratic fields ($\mathbb Q(i)=\mathbb Q(\zeta_4)$ and $\mathbb Q(\sqrt{-3})=\mathbb Q(\zeta_6)$, where $\zeta_n$ is a primitive $n$-th root of unity). Also, as already mentioned, over each of these fields, torsion subgroups appear that appear over no other fields. Note that the rings of integers of both these fields are unique factorization domains.\\

The main results of \cite{fn} are given in the following theorem.
\newtheorem{tm}{Theorem}
\begin{tm}
\begin{itemize}
\item[(i)] Let $E$ be an elliptic curve with rational coefficients. Then $E(\mathbb Q(i))_{tors}$ is either one of the groups from Mazur's Theorem or $\mathbb Z_4 \oplus \mathbb Z_4$.
\item[(ii)] Let $E$ be an elliptic curve defined over $\mathbb Q(i)$. Then $E(\mathbb Q(i))_{tors}$ is either one of the groups from Mazur's Theorem, $\mathbb Z_4 \oplus \mathbb Z_4$ or $\mathbb Z_{13}$.
\item[(iii)] Let $E$ be an elliptic curve with rational coefficients. Then $E(\mathbb Q(\sqrt{-3}))_{tors}$ is either one of the groups from Mazur's Theorem, $\mathbb Z_3 \oplus \mathbb Z_3$ or $\mathbb Z_3 \oplus \mathbb Z_6$.
\item[(iv)] Let $E$ be an elliptic curve defined over $\mathbb Q(\sqrt{-3})$. Then $E(\mathbb Q(\sqrt{-3}))_{tors}$ is either one of the groups from Mazur's Theorem, $\mathbb Z_3 \oplus \mathbb Z_3$, $\mathbb Z_3 \oplus \mathbb Z_6$, $\mathbb Z_{13}$ or $\mathbb Z_{18}$.
\end{itemize}
\end{tm}
While it is not hard to show that (i) and (iii) are best possible, we conjectured that the possible torsion $\mathbb Z_{13}$ could be removed from (ii) and the possible torsions $\mathbb Z_{13}$ and $\mathbb Z_{18}$ could be removed from (iv). In this note we prove this conjecture, thus proving the following theorem.

\newtheorem{tm2}[tm]{Theorem}
\begin{tm2}
\begin{itemize}
\item[(i)] Let $E$ be an elliptic curve defined over $\mathbb Q(i)$. Then $E(\mathbb Q(i))_{tors}$ is either one of the groups from Mazur's Theorem or $\mathbb Z_4 \oplus \mathbb Z_4$.
\item[(ii)] Let $E$ be an elliptic curve defined over $\mathbb Q(\sqrt{-3})$. Then $E(\mathbb Q(\sqrt{-3}))_{tors}$ is either one of the groups from Mazur's Theorem, $\mathbb Z_3 \oplus \mathbb Z_3$ or $\mathbb Z_3 \oplus \mathbb Z_6$.

\end{itemize}
\end{tm2}
\vspace{0.5cm}
\textbf{2 Torsion over $\mathbb Q(i)$ cannot be $\mathbb Z_{13}$}\\
\\
As seen from \cite{ra}, case 2.5.2, page 30, elliptic curves with torsion $\mathbb Z_{13}$ over $\mathbb Q(i)$ exist if and only if the curve $C_1$ defined by
\begin{equation}
C_1: y^2=x^6-2x^5+x^4-2x^3+6x^2-4x+1
\label{he1}
\end{equation}
has points over $\mathbb Q(i)$ satisfying
\begin{equation}
x(x-1)(x^3-4x^2+x+1)\neq 0.
\label{ne1}
\end{equation}

As the points of a hyperelliptic curve have no structure, it is useful to examine the Jacobian variety of a curve (see \cite{ehcc}). Let $J_1$ be the Jacobian of $C_1$, $C_1^{(d)}$ the quadratic twist of $C_1$ by $d$ and $J_1^{(d)}$ the Jacobian of $C_1^{(d)}$.
\newtheorem{tm3}[tm]{Lemma}
\begin{tm3}
$J_1(\mathbb Q(i))\simeq \mathbb Z_{19}.$
\end{tm3}
\emph{Proof:}\\
In MAGMA (\cite{mag}) there exists a implementation of $2$-descent on Jacobians (the {\sf RankBounds} function, see \cite{st} for the algorithm), but unfortunately only over $\mathbb Q$. However,
we can compute $rank(J_1(\mathbb Q))=0$ and $rank(J_1^{(-1)}(\mathbb Q))=0$, and thus
$$ rank(J_1(\mathbb Q(i)))=rank(J_1(\mathbb Q)) + rank(J_1^{(-1)}(\mathbb Q))=0.$$

Again for the torsion of the Jacobian, MAGMA has implemented functions only over the rationals. For an algortihm for computing the torsion, see \cite{po}. The discriminant of the $C_1$ is  $2^{20}\cdot 13^2$, so $2$ and $13$ are the only rational primes with bad reduction.\\
We compute $J_1(\mathbb Q)_{tors}\simeq \mathbb Z_{19}$. If $p$ is a Gaussian prime of good reduction, then the prime-to-$p$ part of  $J_1(Q(i))_{tors}$ injects to $J_1(F_p(i))$. As $p\equiv 3 \pmod 4$ remains prime in $\mathbb Z[i]$ and $F_p(i)\simeq F_{p^2}$, by computing 
$$|J_1(F_{121})|=|J_1(F_{11}(i))|=7^2\cdot 19^2,$$
$$|J_1(F_{529})|=|J_1(F_{23}(i))|=3\cdot 7 \cdot 19 \cdot 229,$$
$$|J_1(F_{961})|=|J_1(F_{31}(i))|=2^8\cdot 3^2 \cdot 19^2,$$
we conclude that $|J_1(Q(i))_{tors}|=19$ and thus $J_1(\mathbb Q(i))\simeq \mathbb Z_{19}$. \qed

\newtheorem{tm4}[tm]{Lemma}
\begin{tm4}
$C_1(\mathbb Q(i))=\{\infty_+, \infty_-, (1, \pm 1), (0, \pm 1)\}.$
\end{tm4}
\emph{Proof:}\\
The 18 non-zero elements of $J_1(\mathbb Q)=J_1(\mathbb Q(i))$, in the Mumford representation that MAGMA uses (see \cite{mum} for details, and \cite{mag2} for the particular implementation in MAGMA) are 
$$(x^2, -2x + 1, 2), (x - 1, -x^3, 2), (x, x^3 - 1, 2), (x, x^3 + 1, 2),(x^2 - x, -1, 2),$$
$$(x^2 - x, -2x + 1, 2), (1, x^3 - x^2, 2), (x^2 - 2x + 1, x, 2),
(x - 1, x^3 - 2, 2),$$ 
$$(x - 1, -x^3 + 2, 2), (x^2 - 2x + 1, -x, 2), (1, -x^3 + x^2, 2), (x^2 - x, 2x - 1, 2),$$
$$(x^2 - x, 1, 2), (x, -x^3 - 1, 2), (x, -x^3 + 1, 2), (x - 1, x^3, 2), (x^2, 2x - 1, 2),$$
and we easily see that the points mentioned are the only ones on this curve. \qed\\

By proving Lemma 4, we have actually proven Theorem 2, (i).\\
\vspace{0.3cm}\\
\textbf{3 Torsion over $\mathbb Q(\sqrt{-3})$ cannot be $\mathbb Z_{13}$ or $\mathbb Z_{18}$}\\
\vspace{0.3cm}\\
To prove that the torsion can not be $\mathbb Z_{13}$, we again have to prove that there are no points on $C_1$ over $\mathbb Q(\sqrt{-3} )$ satisfying (\ref{ne1}).

\newtheorem{tm5}[tm]{Lemma}
\begin{tm5}
$J_1(\mathbb Q(\sqrt{-3}))\simeq \mathbb Z_{19}.$
\end{tm5}
\emph{Proof:\\}
We compute $rank(J_1^{(-3)}(\mathbb Q))=0$, and thus
$$ rank(J_1(\mathbb Q(\sqrt{-3})))=rank(J_1(\mathbb Q)) + rank(J_1^{(-3)}(\mathbb Q))=0.$$
As $5$ and $17$ are rational primes of good reduction that remain prime in $\mathbb Z[\frac{1+\sqrt{-3}}{2}]$, and 
$$|J_1(F_{25})|=|J_1(F_{5}(\sqrt{-3}))|= 19^2,$$
$$|J_1(F_{289})|=|J_1(F_{17}(\sqrt{-3}))|=2^6\cdot 3^2 \cdot 7 \cdot 19,$$
we conclude that $J_1(\mathbb Q(\sqrt{-3}))\simeq \mathbb Z_{19}.$ \qed\\

\newtheorem{tm6}[tm]{Lemma}
\begin{tm6}
$C_1(\mathbb Q(\sqrt {-3}))=\{\infty_+, \infty_-, (1, \pm 1), (0, \pm 1)\}.$
\end{tm6}
\emph{Proof:}\\
The proof is the same as the proof of Lemma 4. \qed\\

As seen from \cite{ra}, case 2.5.6, pages 38--39, elliptic curves with torsion $\mathbb Z_{18}$ over $\mathbb Q(\sqrt {-3})$ exist if and only if the curve $C_2$ defined by
\begin{equation}
y^2=x^6+2x^5+5x^4+10x^3+10x^2+4x+1
\label{he2}
\end{equation}
has points over $\mathbb Q(\sqrt {-3})$ satisfying
\begin{equation}x(x+1)(x^2+x+1)(x^3-3x-1)\neq 0.
\label{ne2}
\end{equation}

Let $J_2$ be the Jacobian of $C_2$, $C_2^{(d)}$ the quadratic twist of $C_2$ by $d$ and $J_2^{(d)}$ the Jacobian of $C_2^{(d)}$.
\newtheorem{tm7}[tm]{Lemma}
\begin{tm7}
$J_2(\mathbb Q(\sqrt{-3}))\simeq \mathbb Z_3 \oplus \mathbb Z_{21}.$
\end{tm7}
\emph{Proof:\\}
In MAGMA (\cite{mag}) we compute $rank(J_2(\mathbb Q))=0$ and $rank(J_2^{(-3)}(\mathbb Q))=0$, and thus
$$ rank(J_2(\mathbb Q(\sqrt{-3})))=rank(J_2(\mathbb Q)) + rank(J_2^{(-3)}(\mathbb Q))=0.$$
The discriminant of $C_2$ is $2^{23} \cdot 3^4$, so $2$ and $3$ are the only rational primes with bad reduction. We find the points in Mumford representation $(x^2 + x + 1, -\sqrt{-3}x - \sqrt{-3}, 2)$ and $(1, -x^3 - x^2, 2)$ on the Jacobian generate a group isomorphic to $\mathbb Z_3 \oplus \mathbb Z_{21}$.

As the rational primes $5$ and $11$ remain prime in $\mathbb Z[\frac{1+\sqrt{-3}}{2}]$, we have 
$$|J_2(F_{25})|=|J_2(F_{5}(\sqrt{-3}))|= 21^2,$$
$$|J_2(F_{121})|=|J_2(F_{11}(\sqrt{-3}))|=2^4\cdot 3^2 \cdot 7 \cdot 13.$$
We conclude that $J_2(\mathbb Q(\sqrt{-3}))_{tors}\simeq \mathbb Z_3 \oplus \mathbb Z_{21}$. \qed

\newtheorem{tm8}[tm]{Lemma}
\begin{tm8}
$C_2(\mathbb Q(\sqrt {-3}))=\{\infty_+, \infty_-, (-1, \pm 1), (0, \pm 1),
(\frac{-1-\sqrt{-3}}{2},
 \pm \frac{-3-\sqrt{-3}}{2} ),$
 $ (\frac{-1+\sqrt{-3}}{2}, \pm \frac{-3+\sqrt{-3}}{2} )\}.$
\end{tm8}
\emph{Proof:}\\
Searching throught the 62 non-zero elements of $J_2(\mathbb Q(\sqrt{-3}))$,
we see that the only points besides the obvious rational ones, satisfy $x^2+x+1=0=y-x+1$, and hence we obtain our result. \qed\\

By Lemma 6 the torsion of an elliptic curve over $\mathbb Q(\sqrt{-3})$ cannot be $\mathbb Z_{13}$ and by Lemma 8, the torsion cannot be $\mathbb Z_{18}$. Thus we have proved Theorem 2.\\

\textbf{Remark}
The fact that $J_1(\mathbb Q)_{tors}\simeq \mathbb Z_{19}$ and $J_2(\mathbb Q)_{tors}\simeq \mathbb Z_{21}$ was proven already 35 years ago by Ogg in \cite{ogg}, and the existence of rational 19-torsion points on $J_1(\mathbb Q)$ was used by Mazur and Tate \cite{mt} to prove the non-existence of rational points on elliptic curves of order 13 over $\mathbb Q$.

\small{FILIP NAJMAN}\\
\small{DEPARTMENT OF MATHEMTICS,\\ UNIVERSITY OF ZAGREB,\\ BIJENI\v CKA CESTA 30, 10000 ZAGREB,\\ CROATIA}\\
\emph{E-mail address:} fnajman@math.hr
\end{document}